\documentclass[12pt]{article}

\usepackage[
    left=1in,
    right=1in,
    top=1in,
    bottom=1in
]{geometry}

\newcommand{\fullFigWidth}{6.5in}
\newcommand{\smallFigWidth}{3.5in}

\newcommand{\reals}{\mathbf{R}}

\newcommand{\diag}{\mathbf{diag}}
\newcommand{\ones}{\mathbf{1}}

\newcommand{\indicator}{\mathbf{I}}
\newcommand{\argmin}{\mathrm{argmin}}

\newcommand{\dom}{\mathbf{dom}}

\newcommand{\prox}{\mathbf{prox}}

\newcommand{\iter}{{(i)}}
\newcommand{\iterhalf}{{(i+1),+}}
\newcommand{\iterplus}{{(i+1)}}
\newcommand{\iterminus}{{(i-1)}}

\usepackage{natbib}
\usepackage{graphicx}
\usepackage{mathtools} 
\usepackage{xcolor} 
\usepackage{verbatim}
\usepackage{hyperref}
\usepackage{amsfonts}

\usepackage{tikz} 
\usetikzlibrary{shapes,positioning,fit}
\tikzset{ell/.style={ellipse,draw,minimum height=0.65cm,minimum width=1cm,inner sep=0.25cm}}

\begin{document}
\title{GPU Accelerated Security Constrained Optimal Power Flow}
\author{
    Anthony Degleris \and
    Abbas El Gamal \and
    Ram Rajagopal
}

\maketitle

\begin{abstract}
We propose a GPU accelerated proximal message passing algorithm for solving contingency-constrained DC optimal power flow problems (OPF). We consider a highly general formulation of OPF that uses a sparse device-node model and supports a broad range of devices and constraints, e.g., energy storage and ramping limits. Our algorithm is a variant of the alternating direction method multipliers (ADMM) that does not require solving any linear systems and only consists of sparse incidence matrix multiplies and vectorized scalar operations. We develop a pure PyTorch implementation of our algorithm that runs entirely on the GPU. The implementation is also end-to-end differentiable, i.e., all updates are automatic differentiation compatible. We demonstrate the performance of our method using test cases of varying network sizes and time horizons. Relative to a CPU-based commercial optimizer, our implementation achieves well over 100$\times$ speedups on large test cases, solving problems with over 500 million variables in under a minute on a single GPU.

\end{abstract}

\section{Introduction}

\textit{Optimal power flow} (OPF) is central to many key problems in power systems operations, economics, and planning~\citep{frank2016introduction}.
To ensure systems are robust to emergency outages, grid operators generally solve \textit{contingency-constrained} OPF problems \citep{capitanescu2011state, frank2016introduction}, which guarantee reliable operation under many potential failure scenarios.
For larger systems, which may interconnect tens of thousands of nodes, generators, and transmission lines, contingency-constrained OPF problems can be extremely large and slow to solve.
Furthermore, increasing renewable penetration has introduced new sources of uncertainty and a need to consider longer time horizons, further increasing the size and complexity of OPF problems.

One way to scale optimization algorithms to larger OPF problems is to leverage vector parallelism using graphics processing units (GPUs).
Unfortunately, this is not straightforward.
Most optimizers---such as interior point methods and simplex methods---require repeatedly solving linear systems, usually via a direct factorization method~\citep{nocedal1999numerical}.
Sparse direct factorizations are naturally implemented sequentially, and existing GPU libraries for these methods do not yet offer significant performance benefits over conventional CPU-based methods~\citep{swirydowicz_linear_2022,shin_accelerating_2024}.

In this work, we use the \textit{alternating direction method of multipliers} (ADMM) to derive a \textit{proximal message passing} algorithm for contingency-constrained OPF problems that does not require solving any linear systems.
We implement this algorithm in PyTorch~\citep{ansel2024pytorch}, a Python-based GPU computing language, and demonstrate that it significantly outperforms commercial CPU-based solvers for large problem instances.

\subsection{Related Work}

A number of works consider using GPUs to accelerate \textit{power flow} calculations, which determine the state of the system without optimizing decisions.
This has been done by applying batch analysis~\citep{araujo_simultaneous_2019, wang_fast_2021}, using iterative linear solvers~\citep{garcia_parallel_2010, li_gpu-based_2017}, or solving the sparse linear system on the GPU directly~\citep{guo_performance_2012, huang_performance_2017, shawlin_gpu-based_2022}.
For \textit{optimal} power flow, GPUs have similarly been used to accelerate interior point methods by evaluating and solving the KKT system on the GPU~\citep{rakai_gpu-accelerated_2014, swirydowicz_gpu-resident_2024} and solving batches of OPF problems in parallel~\citep{kim_leveraging_2021}.
One line of work~\citep{pacaud_feasible_2022, pacaud_condensed-space_2024, shin_accelerating_2024} uses reduced-space methods to convert large sparse problems to smaller dense ones in order to better leverage GPU parallelism.

ADMM is a decomposition technique for splitting complex problems into a series of simpler ones~\citep{boyd2011distributed}.
Various forms of ADMM have been applied to both DC formulations~\citep{kraning2013dynamic,chakrabarti_security_2014,wang_fully-decentralized_2017,yang_parallel_2018,abraham2018admm,javadi_implementation_2019} and AC formulations~\citep{chung_multi-area_2011, sun_fully_2013, erseghe_distributed_2014, erseghe_distributed_2015, shin_hierarchical_2019,sun_two-level_2021, gholami_admm-based_2023} of OPF.
At least one work~\citep{kim_accelerated_2023} considers using GPUs to accelerate ADMM for OPF;
in particular, they combine distributed component-based decompositions with a two-level framework~\cite{sun_two-level_2021} to derive an GPU-amenable ADMM algorithm for single-period AC OPF problems.
Like~\citet{kraning2013dynamic, mhanna_adaptive_2019} and~\citet{kim_accelerated_2023}, our work uses a distributed component-based decomposition, but is unique in considering multi-period, contingency-constrained DC OPF on networks with generic devices.
ADMM is part of a larger class of distributed algorithms that can be used to solve OPF problems; see~\citet{molzahn2017survey} for a general overview.

\subsection{Contribution}

This work has two main contributions.

\paragraph{Message passing for contingency-constrained OPF.}
We extend the model and algorithm from~\citet{kraning2013dynamic} to solve DC linearized optimal power flow problems with contingency-constraints.
We focus on DC linearized OPF because it is convex and widely used in practice~\citep{frank2016introduction}.
This extension preserves the generality of the original model, allowing for devices with a diverse array of costs and constraints.
The resulting algorithm, proximal message passing, does not solve any linear systems and only involves sparse incidence matrix multiples and vectorized scalar operations.

\paragraph{Efficient GPU implementation.}
We show how to efficiently implement this algorithm on a GPU using existing kernels and vectorized proximal operators.
Our open-source implementation\footnote{Available at \url{https://github.com/degleris1/zap}.} is written in pure PyTorch and is end-to-end differentiable, i.e., PyTorch can automatically differentiate the entire algorithm and compute the sensitivity of the solution with respect to problem data.
This GPU implementation greatly outperforms off-the-shelf solvers on realistic grid dispatch problems for low to medium accuracies, e.g., 1\% to 0.01\% suboptimal.
For many applications, such as stochastic production cost modeling or even day-ahead dispatch, this inaccuracy is smaller than the actual uncertainty in the problem data and the error from various modeling assumptions (e.g., the DC power flow approximation).

\section{Optimal Power Flow}
\label{sec:model}

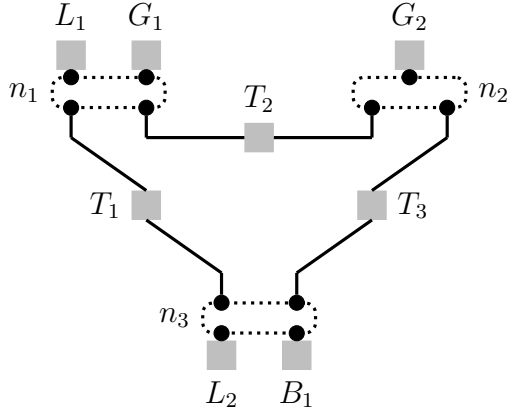
\begin{figure}
\centering
\begin{tikzpicture}[
terminal/.style={circle, draw=black, fill=black, minimum width=0.2cm, inner sep=0pt},
device/.style={rectangle, draw=lightgray, fill=lightgray, very thick, minimum height=0.35cm, minimum width=0.35cm},
node/.style={rectangle, draw=black, dotted, rounded corners, very thick, minimum width=1.5cm, minimum height=0.4cm},
]

\node[device, anchor=south, label=above:$L_1$] (L1) at (-0.5, 0.3) {};
\node[terminal, anchor=north] at (-0.5, 0.31) {};
\node[device, anchor=south, label=above:$G_1$] (G1) at (0.5, 0.3) {};
\node[terminal, anchor=north] at (0.5, 0.31) {};

\node[device, anchor=south, label=above:$G_2$] (G2) at (4.0, 0.3) {};
\node[terminal, anchor=north] at (4.0, 0.31) {};

\node[device, anchor=north, label=below:$L_2$] (L2) at (1.5, -3.3) {};
\node[terminal, anchor=south] at (1.5, -3.31) {};
\node[device, anchor=north, label=below:$B_1$] (B1) at (2.5, -3.3) {};
\node[terminal, anchor=south] at (2.5, -3.31) {};

\node[device, anchor=center, label=left:$T_1$] (T1) at (0.5, -1.5) {};
\node[terminal, anchor=south] (t11) at (-0.5, -0.31) {};
\node[terminal, anchor=south] (t12) at (0.5, -0.31) {};

\node[device, anchor=center, label=above:$T_2$] (T2) at (2.0, -0.6) {};
\node[terminal, anchor=south] (t21) at (3.5, -0.31) {};
\node[terminal, anchor=south] (t22) at (4.5, -0.31) {};

\node[terminal, anchor=north] (t31) at (1.5, -2.69) {};
\node[terminal, anchor=north] (t32) at (2.5, -2.69) {};
\node[device, anchor=center, label=right:$T_3$] (T3) at (3.5, -1.5) {};

\draw[very thick] 
    (t11) edge (-0.5, -0.6) 
    (-0.5, -0.6) edge (T1.north)
    (T1.south) edge (1.5, -2.4) 
    (1.5, -2.4) edge (t31);
\draw[very thick] 
    (t12) edge (0.5, -0.6) 
    (0.5, -0.6) edge (T2.west)
    (T2.east) edge (3.5, -0.6) 
    (3.5, -0.6) edge (t21);
\draw[very thick]
    (t22) edge (4.5, -0.6) 
    (4.5, -0.6) edge (T3.north)
    (T3.south) edge (2.5, -2.4) 
    (2.5, -2.4) edge (t32);

\node[node, label=left:$n_1$] (n1) at (0, 0) {};
\node[node, label=right:$n_2$] (n2) at (4, 0) {};
\node[node, label=left:$n_3$] (n3) at (2, -3) {};
    
\end{tikzpicture}
\caption{
    A 3-bus system with 3 lines, 2 loads, 2 generators, and 1 battery represented in the standard device-node model.
    The system has $J = 11$ terminals (solid black dots), $N = 3$ nodes (dotted boxes), and $D = 8$ devices (gray squares).
    In this example, the devices are $\{ L_1, L_2, G_1, G_2, B_1, T_1, T_2, T_3\}$, and the nodes are $\{ n_1, n_2, n_3\}$ (we use capital letters for the devices to distinguish different types).
    Each device or node is a subset of the terminals, e.g., $n_1 = \{1, 2, 3, 4\}$ and $L_1 = \{ 1\}$, $G_1 = \{ 2 \}$, $T_1 = \{3, 5\}$, $T_2 = \{4, 9\}$.
    Based on Figure~2.1 in~\citet{kraning2013dynamic}.
}
\label{fig:node-diagram}
\end{figure}

In this work, we use the device-node model of electricity networks introduced in~\citet{kraning2013dynamic}.
Consider a network consisting of \textit{terminals} $j \in \mathcal J = \{ 1, \ldots, J \}$ that connect devices $d \in \{ d_1, \ldots, d_D\}$ and nodes $n \in \{n_1, \ldots, n_N\}$, where each device $d$ and node $n$ are subsets of $\mathcal J$.
Each terminal is connected to exactly one device and one node, i.e., the devices and nodes each partition the set of terminals $\mathcal J$.
Devices include, for example, generators, loads, and transmission lines, and generally connect to one, two, or three terminals.
Nodes are electrical buses and may connect to many terminals.
We write $j \in d$ and $j \in n$ if terminal $j$ connects to device $d$ and node $n$, respectively, and let $|d|$ and $|n|$ refer to the number of terminals connected to device $d$ and node $n$, respectively.
We give an example of a 3-bus network with 2 loads, 2 generators, 3 lines, and 1 battery represented in the device-node model in Figure~\ref{fig:node-diagram}.

We consider the operation of the network over $T > 0$ time periods.
Specifically, each terminal $j$ has a time-varying power injection $p_j \in \reals^T$ and phase (angle) $\theta_j \in \reals^T$.
We adopt the convention that $(p_j)_t > 0$ is power flowing from a device to a node at time $t$, i.e., the device is producing power.

\paragraph{Notation.}
We use set-valued indices like $p_n \in \reals^{|n| \times T}$ and $p_d \in \reals^{|d| \times T}$ to denote the power flows of all the terminals connected to a node $n$ and device $d$, respectively, and likewise with phases.
To avoid confusion, we always use $j$ to refer to terminal, $d$ to refer to a device, and $n$ to refer to a node.
Refer to Appendix~\ref{apdx:notation} for a summary of the notation used throughout the text.

\paragraph{}
We define the \textit{average operator} for a terminal $j \in n$ as
\begin{equation}
\label{eq:average-op}
    \bar p_j := \frac{1}{|n|} \sum_{j' \in n} p_{j'}.
\end{equation}
Note that averages are the same for any two terminals connected to the same node: if $j, j' \in n$ then $\bar p_j = \bar p_{j'}$.
Like $p$, we write $\bar p_n \in \reals^{|n| \times T}$ and $\bar p_d \in \reals^{|n| \times T}$ to refer to the average power injections for a group of terminals connected to a node or device.
By definition, $\bar p_n$ consists of $|n|$ identical copies of the average power injection at node $n$.
Each node $n$ in the network is constrained to satisfy \textit{power balance}, $\bar p_n = 0$, i.e., that the average power at node $n$ must be zero.
We also define the \textit{residual operator} as
\begin{equation}
\label{eq:residual-op}
    \tilde \theta_j := \theta_j - \bar \theta_j.
\end{equation}
Nodes must satisfy \textit{phase consistency}, $\tilde \theta_n = 0$, which says that all phases at node $n$ are equal.

Finally, each device $d$ has a closed proper convex cost function $f_d : \reals^{|d| \times T} \times \reals^{|d| \times T} \rightarrow \reals \cup \{ \infty \}$ that specifies both the costs and constraints of operating that device with a given power and phase schedule.
The functions $f_d$ may take on infinite values to encode constraints, and $\dom f_d = \{ (p_d, \theta_d) \mid f_d(p_d, \theta_d) < \infty \}$ is the feasible set of schedules for a device.
The \textit{DC-linearized optimal power flow (OPF)} problem is
\begin{equation}
\label{eq:dispatch}
\begin{array}{ll}
    \textrm{minimize} \quad
    & \sum_{d} f_d(p_d, \theta_d) \\[0.5em]
    \textrm{subject to} \quad
    & \bar p = 0, \quad \tilde \theta = 0,
\end{array}
\end{equation}
where the variables are $p = (p_1, \ldots, p_J) \in \reals^{J \times T}$ and $\theta = (\theta_1, \ldots, \theta_J) \in \reals^{J \times T}$.

\subsection{Device Examples}

Below, we give examples of several common device cost functions.

\paragraph{Generator.}
A generator is a single-terminal device, $|d| = 1$, with simple box constraints $p^{\min} \leq p_1 \leq p^{\max}$, where $p^{\min}, p^{\max} \in \reals^T$.
The generator has quadratic operation cost $a p_1^2 + b p_1$.
Its device cost function is
\begin{equation*}
    f^{\textrm{gen}}(p_1, \theta_1) = a p_1^2 + b p_1 + \indicator \{ p^{\min} \leq p_1 \leq p^{\max} \},
\end{equation*}
where $\indicator \{ A x \leq b \}$ is the convex indicator function on the constraint $Ax \leq b$, which evaluates to $\infty$ if the constraint is not satisfied and zero otherwise.
Curtailable and flexible loads can be modeled as generators with $p^{\min} < p^{\max} \leq 0$.

\paragraph{Fixed load.}
Fixed loads are represented by a simple single-terminal device with equality constraints,
\begin{equation*}
    f^{\textrm{load}}(p_1, \theta_1) = \indicator \{ p_1 = p^{\textrm{load}} \},
\end{equation*}
where $p^{\textrm{load}} \in \reals^T$.

\paragraph{AC transmission line.}
Transmission lines are two-terminal devices that transfer power from one node to another.
Under the \textit{DC power flow approximation}, a (lossless) AC transmission line has device cost function,
\begin{equation*}
    f^{\textrm{ac}}(p_1, p_2, \theta_1, \theta_2) = \indicator \{ p_1 + p_2 = 0 \} + \indicator \{ -u \ones \leq p_2 \leq u \ones \} + \indicator \{ p_2 = b (\theta_1 - \theta_2) \},
\end{equation*}
where $u \in \reals$ is the line flow limit and $b \in \reals$ is the line susceptance.

\paragraph{Battery.}
A battery is a complex single-terminal device involving local variables: the charge schedule $c \in \reals^T$, the discharge schedule $\delta \in \reals^T$, and the state of charge $s \in \reals^{T+1}$.
The battery has discharge cost $\alpha \in \reals_+$, a charging efficiency $\beta \in (0, 1]$, a power capacity $P \in \reals_+$, and a duration $\lambda \in \reals_+$.
The device cost function is
\begin{align*}
    f^{\textrm{store}}(p_1, \theta_1)
    = \min_{c, d, s} \Big( &\alpha \ones^T d
    + \indicator\{ p_1 = \delta - c \}
    + \indicator\{ s_{t+1} = s_t + \beta c_t - \delta_t \} \\
    &+ \indicator\{ 0 \leq c \leq P \ones \}
    + \indicator\{ 0 \leq \delta \leq P \ones \}
    + \indicator\{ 0 \leq s \leq \lambda P \ones \} \\
    &+ \indicator\{ s_1 = 0 \} + \indicator\{ s_{T+1} = 0 \} \Big).
\end{align*}
The first constraint is device power balance, and the second constraint describe the state of charge evolution.
The third, fourth, and fifth constraints are limits on charge, discharge, and state of charge.
The final two constraints are initial and final conditions on the battery.

\paragraph{Other devices.}
Many other common electrical devices can be represented under this framework, such as HVDC transmission lines, transformers, generators with ramping limits, and simple models of lossy transmission lines.
Additionally, other energy products, e.g., ancillary services or gas flow, can be modeled by designating special nodes to be carriers of said products.

\paragraph{Non-convex devices.}
Some devices, such as generators with unit commitment constraints or lossy transmission lines, have a non-convex device cost function.
We do not consider these devices in this work.
This can be addressed in a few different ways, which we discuss as possible directions for future work in Section~\ref{sec:conclusion}.

\subsection{Vectorization}
\label{sec:vectorize}


Previous work~\citep{kraning2013dynamic} uses the device-node model to derive fully distributed algorithms for solving grid dispatch problems.
However, if the number of unique device types is large, these distributed algorithms will perform poorly on vector processing architectures such as GPUs.
Therefore, we impose the requirement that each device $d$ is a member of a type set $\ell \in \{\ell_1, \ldots, \ell_L\}$, where $L \ll D$.
Each type set $\ell$ consists of devices $d \in \ell$; the quantity $|\ell|$ refers to the number of devices of type $\ell$.
When $d \in \ell$, we require that the device has $|d| = \tau_\ell$ terminals and device cost function,
\begin{equation*}
    f_d(p_d, \theta_d) = g_{\ell}(p_d, \theta_d; x_d),
\end{equation*}
where $g_\ell : \reals^{\tau_\ell \times T} \times \reals^{\tau_\ell \times T}  \times \reals^{r_\ell} \rightarrow \reals \cup \{ \infty \}$ is the type cost function.
Each device has a \textit{parameter} $x_d \in \reals^{r_\ell}$, which is a fixed value that determines the specific characteristics of its device.
For example, for an AC transmission line, the parameter is $x_d = (u_d, b_d) \in \reals^2$, i.e., the line flow limit and susceptance.

Just as we use $p_n \in \reals^{|n| \times T}$ and $p_d \in \reals^{|d| \times T}$ to index groups of terminals connected to a particular node or device, we denote $p_\ell \in \reals^{|\ell| \times \tau_\ell \times T}$ to all the terminals connected to a device of type $\ell$.
We then define
\begin{equation*}
    G_{\ell}(p_\ell, \theta_\ell; x_\ell)
    = \sum_{d \in \ell} g_{\ell}(p_d, \theta_d; x_d),
\end{equation*}
so the objective of Problem~\eqref{eq:dispatch} can be rewritten as
\begin{equation*}
    \sum_{\ell=1}^L G_{\ell}(p_\ell, \theta_\ell; x_\ell).
\end{equation*}
As we shall see in Section~\ref{sec:implementation}, this allows related computations for all devices of type $\ell$ to be batched and executed in parallel.
In practice, each device type $\ell$ might correspond to one of the examples in Section~\ref{sec:model}, and $L$ will be a small constant, e.g., $L = 5$.

\subsection{Contingencies}

Suppose there are contingency events $k \in \{1, \ldots, K\}$ during which some of the device parameters change.
Specifically, devices of the first type change their parameters from $x_d$ to $x_{kd} \in \reals^{r_\ell}$; we will use $x_{0d}$ to denote their base case parameters when no contingency event occurs.
These devices may also change their power and phase schedules in response to parameter changes.
Devices of the remaining $L - 1$ types have fixed parameters, $x_{kd} = x_d$, and cannot change their power or phase schedules.
The \textit{contingency-constrained optimal power flow problem} is
\begin{equation}
\label{eq:contingency-dispatch}
\begin{array}{ll}
    \textrm{minimize} \quad
    & \frac{1}{K+1} \sum_{k=0}^{K} \left( G_\ell(p_{k \ell_1}, \theta_{k \ell_1}; x_{k \ell_1})
    + 
    \sum_{\ell \neq \ell_1} G_\ell(p_{k \ell}, \theta_{k \ell}; x_{k}) \right)
    \\[0.5em]
    \textrm{subject to} \quad
    & \bar p_k = 0, \quad \tilde \theta_k = 0, \\
    & p_{k \ell} = p_{0 \ell}, \quad \ell = 2, \ldots, L,\quad k = 1, \ldots, K, \\
    & \theta_{k \ell} = \theta_{0 \ell}, \quad \ell = 2, \ldots, L,\quad k = 1, \ldots, K,
\end{array}
\end{equation}
where the variables are $p_0, \ldots, p_K \in \reals^{J \times T}$ and $\theta_0, \ldots, \theta_K \in \reals^{J \times T}$.
Here, we use 
$p_{k \ell} := (p_k)_{\ell} \in \reals^{|\ell| \times \tau_\ell \times T}$ to index the power schedule of 
type $\ell$ for contingency $k$.
Contingency-constrained OPF problems include, as a special case, the \textit{security constrained economic dispatch problem} (SCED), which is widely solved in practice when scheduling day-ahead power operations~\citep{frank2016introduction}.
It is straightforward to generalize our framework here to problems where multiple device types may change their data and respond to contingency events;
we focus on single type contingencies for ease of presentation.

\paragraph{Example: $N-1$ security constraints.}
The canonical example of contingency-constrained OPF is the $N-1$ security-constrained OPF problem, in which solutions must be robust to any single transmission line outage.
In our framework, suppose the first device type is an AC line, $g_1 = f^{\textrm{ac}}$, of which there are $K$ devices.
The base case parameters are the transmission line capacities and susceptances $x_1^0 = (u, b) \in \reals^{2K}$.
The parameters for contingency $k$ are $x_1^k = (u \odot (\ones - e_k), b \odot (\ones - e_k))$, where $e_k \in \reals^K$ is the $k$th standard basis vector and $\odot$ is elementwise multiplication.
In other words, the parameter for contingency $k$ is the same as the base case, except the capacity and susceptance of the line $k$ is now zero.

\section{Message Passing Algorithm}
\label{sec:algorithm}

In this section, we introduce our message passing algorithm for solving contingency-constrained OPF.
Specifically, we reformulate~\eqref{eq:contingency-dispatch} by replacing constraints with indicator functions and introducing additional variables.
We then show that the updates in the alternating direction method of multipliers (ADMM) can be simplified to obtain a proximal \textit{message passing algorithm}.
Refer to~\citet{kraning2013dynamic} for the same algorithm in non-contingency case, and to~\citet{boyd2011distributed} for details on ADMM in general.

\subsection{Problem Reformulation}

We reformulate Problem~\eqref{eq:contingency-dispatch} as
\begin{equation}
\label{eq:admm-problem-full}
\begin{array}{lll}
    \textrm{minimize} \quad
    & \frac{1}{K+1} \sum_{k=0}^{K} \left(
    \indicator\{ \bar z_k = 0\} + \indicator \{ \tilde \xi_k = 0 \}
    + G_\ell(p_{k \ell_1}, \theta_{k \ell_1}; x_{k \ell_1})
    + \sum_{\ell \neq \ell_1} G_\ell(p_{k \ell}, \theta_{k \ell}; x_k) \right)
    \\[1em]
    \textrm{subject to} \quad
    & p_k = z_k, \quad \theta_k = \xi_k, \quad k = 1, \ldots, K \\[0.5em]
    & p_{k \ell} = p_{0 \ell}, \quad \ell = \ell_2, \ldots, \ell_L,\quad k = 1, \ldots, K, \\[0.5em]
    & \theta_{k \ell} = \theta_{0 \ell}, \quad \ell = \ell_2, \ldots, \ell_L, \quad k = 1, \ldots, K,
\end{array}
\end{equation}
where the variables are $p_k, \theta_k \in \reals^{J \times T}$ and $z_k, \xi_k \in \reals^{J \times T}$ for $k = 0, \ldots, K$.
The variables $z_k \in \reals^{J \times T}$ and $\xi_k \in \reals^{J \times T}$ are copies of the power and phase variables, respectively, which are used to separate the node constraint functions from the device cost functions $f_d$.

Since the power and phase schedules do not change for devices of type $\ell \neq \ell_1$ across different contingencies, we can replace $p_{0 \ell}, \ldots, p_{K \ell}$ with a single variable $p_\ell \in \reals^{|d| \times T}$.
This lets us reformulate the constraints of Problem~\eqref{eq:admm-problem-full} as two sets of constraints between $p, \theta$ and their copies,
\begin{equation}
\label{eq:admm-problem}
\begin{array}{lll}
    \textrm{minimize} \quad
    & \sum_{\ell \neq \ell_1} G_\ell(p_\ell, \theta_\ell; x_\ell)
    + \frac{1}{K+1} \sum_{k=0}^{K} \Big( G_1(p_{1 k}, \theta_{1 k}; x_{1k})
    + \indicator \{ \bar z_k = 0 \} + \indicator \{ \tilde \xi_k = 0 \}
    \Big)
    \\[1em]
    \textrm{subject to} \quad
    & (p_{k \ell}, \theta_{k \ell}) = (z_{k \ell}, \xi_{k \ell}),
        \quad \ell = 1, \quad k = 1, \ldots, K, \\[0.5em]
    & (p_\ell, \theta_\ell) = (z_{k \ell}, \xi_{k \ell}),
        \quad \ell = 2, \ldots, L, \quad k = 1, \ldots, K,  \\[0.5em]
\end{array}
\end{equation}
where the variables are the non-contingency power and phase schedules $p_\ell, \theta_\ell \in \reals^{|\ell| \times \tau_\ell \times T}$ for $\ell \neq \ell_1$, the (changeable) contingency power and phase schedules $p_{k \ell}, \theta_{k \ell} \in \reals^{|\ell| \times \tau_1 \times T}$ for $\ell = 1$ and $k = 0, \ldots, K$, and the power and phase copies $z_k, \xi_k \in \reals^{J \times T}$ for $k = 0, \ldots, K$.
The objective of~\eqref{eq:admm-problem} is decomposed into two parts: the non-contingency device costs, summed over all device types $\ell > 2$, and the sum of the contingency device costs plus the constraint indicator functions, summed over all  $K+1$ contingencies.
Therefore, a \textit{single} power and phase schedule for each non-contingency device must support power balance, phase consistency, and contingency device constraints for \textit{all} contingencies.

\paragraph{Augmented Lagrangian.}
Following~\citet{boyd2011distributed} and \citet{kraning2013dynamic}, the augmented Lagrangian of~\eqref{eq:admm-problem} is, up to a constant,
\begin{align*}
    L(p, \theta, z, \xi, u, v)
    = &\sum_{\ell \neq \ell_1} G_\ell(p_\ell, \theta_\ell; x_\ell) \\
    &+ \frac{1}{K+1} \sum_{k=0}^{K} \Big( G_1(p_{k \ell_1}, \theta_{k \ell_1}; x_{k \ell_1})
    + \indicator \{ \bar z_k = 0 \} + \indicator \{ \tilde \xi_k = 0 \}
    \Big) \\
    &+ \frac{\rho_p}{2} \sum_k \Big( \| p_{k 1} - z_{k 1} + u_{k 1} \|_2^2 + \sum_{\ell \neq \ell_1} \| p_\ell - z_{k \ell} + u_{k \ell} \|_2^2 \Big) \\
    &+ \frac{\rho_\theta}{2} \sum_k \Big( \| \theta_{k 1} - \xi_{k 1} + v_{k 1} \|_2^2 + \sum_{\ell \neq \ell_1} \| \theta_\ell - \xi_{k \ell} + v_{k \ell} \|_2^2 \Big),
\end{align*}
where $u_k, v_k \in \reals^{J \times T}$ are the scaled dual variables of the power balance and phase consistency constraints, i.e., $\rho_p u_k$ is the dual variable of $p_k = z_k$ and $\rho_\theta v_k$ is the dual variable of $\theta_k = \xi_k$.

\subsection{Message Passing}

The ADMM updates for~\eqref{eq:admm-problem} on iteration $\iter$ are:
\begin{enumerate}
\item
\textit{Contingency proximal schedule updates.}
For contingencies $k = 0, \ldots, K$ and devices $d \in \ell_1$,
\begin{align*}
    (p_{k d}^\iterplus, \theta_{k d}^\iterplus) := \argmin_{p, \theta} \Big(
        g_\ell(p, \theta; x_{kd})
        + \frac{\rho_p}{2} \| p - z_{kd}^\iter + u_{kd}^\iter \|_2^2
        + \frac{\rho_\theta}{2} \| \theta - \xi_{kd}^\iter + v_{kd}^\iter  \|_2^2
    \Big),
\end{align*}

\item
\textit{Non-contingency proximal schedule updates.}
For types $\ell = \ell_2, \ldots, \ell_L$ and devices $d \in \ell$,
\begin{align*}
    (p_d^\iterplus, \theta_d^\iterplus) := \argmin_{p, \theta} \Big(
        g_\ell(p, \theta; x_d)
        + \sum_k \left( \frac{\rho_p}{2} \| p - z_{kd}^\iter + u_{kd}^\iter \|_2^2
        + \frac{\rho_\theta}{2} \| \theta - \xi_{kd}^\iter + v_{kd}^\iter  \|_2^2 \right)
    \Big),
\end{align*}

\item \textit{Duplicate variable updates.}
\begin{align*}
    z_k^{\iterplus} := \argmin_{z} \left(  \indicator \{ \bar z_k = 0 \} + \frac{\rho_p}{2} \| z - p_k^\iterplus - u_k^\iter \|_2^2 \right), \quad k = 0, \ldots, K, \\
    \xi_k^\iterplus := \argmin_{\xi} \left( \indicator \{ \tilde \xi_k = 0 \} + \frac{\rho_\theta}{2} \| \xi - \theta_k^\iterplus - v_k^\iter \|_2^2 \right), \quad k = 0, \ldots, K,
\end{align*}

\item \textit{Scaled price updates.}
\begin{align*}
    u_k^\iterplus := u_k^\iter + p_k^\iterplus - z_k^\iterplus, \quad k = 0, \ldots, K, \\
    v_k^\iterplus := v_k^\iter + \theta_k^\iterplus - \xi_k^\iterplus, \quad k = 0, \ldots, K.
\end{align*}
\end{enumerate}

We can make several simplifications to these update rules.
First, the sums of squared norms in the non-contingency updates can be simplified to a single term plus a constant, allowing us to rewrite the second update rule as
\begin{align*}
    (p_d^\iterplus, \theta_d^\iterplus) := \argmin_{p, \theta} \Big(
        g_\ell(p, \theta; x_d)
        &+  ((K+1) \rho_p / 2) \| p - \mathbf{avg}_k (z_{k d}^\iter - u_{k d}^\iter) \|_2^2 \\
        &+ ((K+1) \rho_\theta / 2) \| \theta - \mathbf{avg}_k ( \xi_{k d}^\iter - v_{k d}^\iter )  \|_2^2
    \Big),
\end{align*}
where $\mathbf{avg}_k(z_k)$ averages across the $K+1$ contingencies.
Second, the duplicate variable updates are just projections onto the constraints $\{ \bar z = 0 \}$ and $\{ \tilde \theta = 0 \}$ and have closed-form updates
\begin{align*}
    z_k^\iterplus &:= \tilde p_k^\iterplus + \tilde u_k^\iter, \quad k = 0, \ldots, K, \\
    \xi_k^\iterplus &:= \bar \theta_k^\iterplus + \bar v_k^\iter, \quad k = 0, \ldots, K.
\end{align*}
To see this, we simply solve the proximal minimization analytically and use the basic properties of the average and residual operators defined in Equations~\eqref{eq:average-op} and~\ref{eq:residual-op}.
Refer to Appendix~\ref{apdx:average-resid} for further detail.
Finally, we can assume that the dual variables satisfy $\tilde u_k = 0$ and $\bar v_k = 0$ at optimality (see Appendix~\ref{apdx:dual}).
If this also holds for iterate~$\iter$, then the duplicate variable updates on iterate $\iterplus$ are simply $z_k^\iterplus := \tilde p_k^\iterplus$ and $\xi_k^\iterplus := \bar \theta_k^\iterplus$.
Substituting these expressions into the price updates,
\begin{align*}
    u_k^\iterplus 
        &= u_k^\iter + p_k^\iterplus - \tilde p_k^\iterplus
        = u_k + \bar p_k^ \iterplus, \\
    v_k^\iterplus 
        &= v_k^\iter + \theta_k^\iterplus - \bar \theta_k^\iterplus 
        = v_k^\iter + \tilde \theta_k^\iterplus.
\end{align*}
Since the average and residual operators are orthogonal (see Appendix~\ref{apdx:average-resid}), the above expressions imply $\tilde u_k^\iterplus = 0$ and $\bar v_k^\iterplus = 0$.
By induction, this means that if $\tilde u_k^{(0)} = 0$ and $\bar v_k^{(0)} = 0$, then $\tilde u_k^\iter = 0$ and $\bar v_k^\iter = 0$ for every iterate thereafter.
Substituting the expressions for $z_k^\iterplus$ and $\xi_k^\iterplus$ into the updates for $(p, \theta)$ as well, we obtain the following \textit{message passing algorithm}:

\begin{enumerate}
\item
\textit{Contingency proximal schedule updates.}
For $k = 0, \ldots, K$ and $d \in \ell_1$,
\begin{align*}
    (p_{k d}^\iterplus, \theta_{k d}^\iterplus) := \prox_{g_\ell, \rho_p, \rho_\theta}\left(\tilde p_{k d}^\iter - u_{k d}^\iter, \bar \theta_{k d}^\iter - v_{k d}^\iter  \right).
\end{align*}

\item
\textit{Non-contingency proximal schedule updates.}
For $\ell = \ell_2, \ldots, \ell_L$ and $d \in \ell$,
\begin{align*}
    (p_d^\iterplus, \theta_d^\iterplus) := \prox_{g_\ell, \rho_p (K + 1), \rho_\theta (K + 1)}\left(
        p_d^\iter - \mathbf{avg}_k( \bar p_{k d}^\iter - u_{k d}^\iter ),\
        \mathbf{avg}_k( \bar \theta_{k d}^\iter - v_{k d}^\iter )
    \right).
\end{align*}

\item \textit{Scaled price updates.}
\begin{align*}
    u_k^\iterplus := u_k^\iter + \bar p_k^\iterplus, \quad k = 0, \ldots, K, \\
    v_k^\iterplus := v_k^\iter + \tilde \theta_k^\iterplus, \quad k = 0, \ldots, K.
\end{align*}
\end{enumerate}

Above, the proximal operator for a function $g$ and scalars $\rho_p, \rho_\theta$ is defined as
\begin{equation*}
    \prox_{g, \rho_1, \rho_2}(x, y) = \argmin_{p, \theta} \left( g(p, \theta) + (\rho_1 / 2) \| p - x \|_2^2 + (\rho_2 / 2) \| \theta - y \|_2^2 \right).
\end{equation*}
The quantities $\rho_p$ and $\rho_\theta$ are the \textit{penalty parameters} that define how quickly that power and phase schedules may change between iterates.
In particular, the larger $\rho_p$ is, the more the power schedule will be penalized for deviating from the ``target schedule'' passed to the proximal operator; likewise for the $\rho_\theta$ and the phase schedule.
The choice of $\rho_p$ and $\rho_\theta$ will depend on the specific characteristics and scaling of the data.

\subsection{Convergence}
\label{sec:convergence}

The primal and dual residuals for contingency $k$ are,
\begin{align*}
    r_{\textrm{primal}, k}^\iter
    &= (\bar p_k^\iter, \tilde \theta_k^\iter), \\
    r_{\textrm{dual}, k}^\iter
    &= (
        \rho_p (\tilde p_k^\iter - \tilde p_k^\iterminus),
        \rho_\theta (\bar \theta_K^\iter - \bar \theta_K^\iterminus)
    ),
\end{align*}
and the total residuals are $r_{\textrm{primal}} = (r_{ \textrm{primal}, 0}, \ldots, r_{\textrm{primal}, K})$ and $r_{\textrm{dual}} = (r_{ \textrm{dual}, 0}, \ldots, r_{\textrm{dual}, K})$.
The primal residual is simply the power imbalance and phase inconsistency across all contingencies, i.e., the total primal infeasibility.
The dual residual is the change between iterates in the power and phase schedules projected onto the power balance and phase consistency constraints;
a small change indicates dual convergence and an optimal solution.
We terminate the algorithm when the root mean square residuals are less than some tolerance,
\begin{equation}
\label{eq:converge}
    \max\left( \| r_{\textrm{primal}} \|_2, \| r_{\textrm{dual}} \|_2 \right) \leq \epsilon \sqrt{2 J T (K+1)},
\end{equation}
where $\epsilon > 0$ is the absolute tolerance that will depend on the specific problem data, e.g., the total load, the units for cost, and so on.
Because our problem is convex, message passing is guaranteed to converge if a feasible solution exists;
see~\citet{boyd2011distributed} for further detail on the theoretical convergence properties of ADMM.

\subsection{Warm Starts}
\label{sec:warm-start}

By default, we initialize the algorithm with $(p^{(0)}, \theta^{(0)}, u^{(0)}, v^{(0)}) = 0$.
However, any point $(p, \theta, u, v)$ such that $\tilde u = 0$ and $\bar v = 0$ is a valid initialization.
In particular, we can initialize the algorithm using the solution to a similar problem.
For example, if we are solving one day of hourly operations (so $T = 24$), then we can initialize the algorithm with the solution to the previous 24 hours of data.
In general, warm starts are useful both for operations problems, where the solution from the previous interval is readily available, and planning problems, where many similar simulations are repeatedly solved with different network parameters.
If the solutions are similar, this should significantly reduce the number of iterations required to converge to a particular accuracy.

\subsection{Adaptive Step Size}
\label{sec:adaptive-penalty}

Choosing the penalty parameters $\rho_p$ and $\rho_\theta$ can be a priori difficult since they depend on the choice of device cost functions $g_\ell$ and their data $x_\ell$.
Moreover, the optimal choice of step size can change over the course of the algorithm.
We use a simple adaptive penalty update from~\citet{boyd2011distributed},
\begin{align*}
    \rho_{x}^{\iterplus} := \begin{cases}
        \gamma \rho_x^\iter
            & \| r_{\textrm{primal}, x}^\iter \|_2 > \mu \| r_{\textrm{dual}, x}^\iter \|_2 \\
        \rho_x^\iter / \gamma
            & \| r_{\textrm{dual}, x}^\iter \|_2 > \mu \| r_{\textrm{primal}, x}^\iter \|_2 \\
        \rho_x^\iter
            & \textrm{otherwise}
    \end{cases}
\end{align*}
where $x \in \{p, \theta \}$ and
\begin{align*}
    r_{\textrm{primal}, p}^\iter &= (\bar p_0^\iter, \ldots, \bar p_K^\iter), \\
    r_{\textrm{primal}, \theta}^\iter &= (\tilde \theta_0^\iter, \ldots, \tilde \theta_K^\iter), \\
    r_{\textrm{dual}, p}^\iter &= \rho_p (\tilde p_0^\iter - \tilde p_0^\iterminus, \ldots, \tilde p_K^\iter - \tilde p_K^\iterminus), \\
    r_{\textrm{dual}, \theta}^\iter &= \rho_\theta (\bar \theta_0^\iter - \bar \theta_0^\iterminus, \ldots, \bar \theta_K^\iter - \bar \theta_K^\iterminus).
\end{align*}
Since $u$ and $v$ are the scaled dual variables, they must be rescaled each time the penalty parameters are updated,
\begin{align*}
    u^\iterplus &:= \left( \rho_p^\iter / \rho_p^\iterplus \right) u^\iterplus, \\
    v^\iterplus &:= \left( \rho_\theta^\iter / \rho_\theta^\iterplus \right) v^\iterplus.
\end{align*}
This update rule adapts the penalty parameters for power and phase separately based on the primal and dual residuals associated with those variables and constraints.
This has the effect of adapting the relative scaling between power and phase variables, i.e., the unit used for line susceptance.
In the language of power systems, changing $\rho_\theta$ for a fixed $\rho_p$ is equivalent to changing the base value for impedance in the per unit system, while leaving the base value for power fixed.

In practice, we find that this adaptive scheme works well with $\mu = 2$ and $\gamma = 1.1$.
However, updating $\rho_p$ or $\rho_\theta$ may mean recomputing cached quantities used to evaluate proximal operators.
Therefore, we only update the penalty parameters every $i_{\textrm{adapt}}$ iterations; we set $i_{\textrm{adapt}} = 10$ in all our experiments.

\subsection{Over-relaxation}

\textit{Over-relaxation} is a technique to accelerate the convergence of ADMM~\citep{boyd2011distributed}.
In the $(z, \xi)$ and $(u, v)$ updates, we replace $(p_k^\iter, \theta_k^\iter)$ with,
\begin{align*}
    p_k^\iterhalf &:= \alpha p_k^\iterplus + (1 - \alpha) z_k^\iter, \\
    \theta_k^\iterhalf &:= \alpha \theta_k^\iterplus + (1 - \alpha) \xi_k^\iter,
\end{align*}
where $\alpha \in [1, 2)$ is called the \textit{relaxation parameter}.
When $\alpha = 1$, we recover standard (unrelaxed) ADMM.
The over-relaxed $(z, \xi)$ updates are
\begin{align*}
    z_k^\iterplus &:= \alpha \tilde p_k^\iterplus + (1 - \alpha) z_k^\iter, \\
    \xi_k^\iterplus &:= \alpha \bar \theta_k^\iterplus + (1 - \alpha) \xi_k^\iter,
\end{align*}
and the scaled price updates are
\begin{align*}
    u_k^\iterplus &:= u_k^\iter + \alpha \bar p_k^\iterplus, \\
    v_k^\iterplus &:= v_k^\iter + \alpha \tilde \theta_k^\iterplus.
\end{align*}
Because we need to keep track of $(z, \xi)$ from the previous iteration, we cannot simplify ADMM to the message passing algorithm; instead, we use the original ADMM updates and keep explicit copies of the duplicate variables.
In practice, we find $\alpha = 1.5$ works well.

\section{GPU Implementation}
\label{sec:implementation}

We implement the entire message passing algorithm on the GPU to take advantage of parallelism in many of the update rules.
There are two types of operations that need to be implemented: proximal operators, evaluated in parallel across many devices, and averaging and residual operations applied across nodes and contingencies.
In this section, we show how to implement each of these operations efficiently on the GPU.
Then, we explain how the solver output can be differentiated with respect to device parameters using an automatic differentiation library, such as those implemented in PyTorch~\citep{ansel2024pytorch}, TensorFlow~\citep{tensorflow2015-whitepaper}, or JAX~\citep{jax2018github}.

\subsection{Data Model}

For each device type $\ell$, we represent its variables as $\tau_\ell$ tensors of size $(K+1) \times |\ell| \times T$; each tensor describes the variable associated with one terminal per device of that type.
For quantities shared across contingencies, such as $p_\ell$ for $\ell \geq 2$, the dimension of the first axis is set to 1.
Quantities that are constant for all terminals in a node, like $\bar p$ or $u$, are stored as tensors of size $(K+1) \times N \times T$.
Finally, device data $x_\ell$ are represented as tensors of size $(K+1) \times |\ell| \times r_\ell$.
The benefit of grouping data by device type is that all devices of type $\ell$ have the same cost function and the same number of terminals.
Therefore, operations for a given device type can be \textit{vectorized}, i.e., batch processed across many processors in parallel.

\subsection{Proximal Operators}

Recall the definition of the proximal operator is
\begin{equation*}
    \prox_{g, \rho_1, \rho_2}(x, y) = \argmin_{p, \theta} \left( g(p, \theta) + (\rho_1 / 2) \| p - x \|_2^2 + (\rho_2 / 2) \| \theta - y \|_2^2 \right).
\end{equation*}
In step~1 of the message passing algorithm, we evaluate a $(K+1) |\ell_1|$ proximal operators for the function $g_1(\cdot, \cdot; x_d) : \reals^T \times \reals^T \rightarrow \reals$.
In step~2 of the message passing algorithm, we evaluate $|\ell|$ proximal operators for functions $g_\ell(\cdot, \cdot; x_d)$ for each $\ell = \ell_2, \ldots, \ell_L$.
In total, evaluating the proximal operators in steps~1 and~2 require solving $(K+1) |\ell_1| + (D - |\ell_1|)$ minimization problems over $\reals^{2T}$.

\paragraph{Analytical updates.}
Many proximal operators have simple closed-form expressions.
For example, a generator has proximal update,
\begin{equation*}
    \prox_{f^{\textrm{gen}}, \rho}(z)
    =
    \mathrm{\Pi}_{[p^{\min}, p^{\max}]}\left( \frac{2a + \rho}{\rho z - b} \right),
\end{equation*}
where $\mathrm{\Pi}_{[a, b]}$ is the elementwise projection onto the box $[a, b]$ and the terms for $\theta$ are omitted since they are not involved in the device cost functions or constraints (i.e., $\theta = y$).
For an AC line, the proximal update has solution $\prox_{f^{\textrm{ac}}, \rho_p, \rho_\theta}(z, \xi) = (p_1, p_2, \theta_1, \theta_2)$ where
\begin{align*}
    p_2 &= \mathrm{\Pi}_{[-u, u]} \left( \frac{\rho_p (z_2 - z_1) + \rho_\theta (\xi_1 - \xi_2) / 2b}{2 (\rho_p + \rho_\theta + 1/4b^2)} \right), \\
    p_1 &= -p_2, \\
    \theta_2 &= (1/2) (\xi_2 + \xi_1) - p_2 / 2 b, \\
    \theta_1 &= \theta_2 + p_2 / b.
\end{align*}

\paragraph{Approximate iterative updates.}
Some devices, such as batteries or generators with ramping constraints, have more complex device cost functions whose proximal operators cannot be evaluated analytically.
Generally, this occurs when the devices have constraints that couple variables across time.
In the most general setting, we assume the device cost $f$ is the minimum of some quadratic program,
\begin{align}
\label{eq:prox-general-form}
    f(p, \theta) = \min_{s} \left( s^T Q s + q^T s + \indicator\{ A_1 p + A_2 \theta + A_3 s \leq b \} \right),
\end{align}
where $s \in \reals^{\mu}$ is a vector of local variables.
Evaluating $\prox_{f, \rho_p, \rho_\theta}(x, y)$ requires solving a quadratic program with linear constraints and $2T + \mu$ variables.
When $T$ and $\mu$ are reasonably small, these quadratic programs can be solved approximately using a few iterations of an iterative algorithm.
In our implementation, we use 10 to 50 iterations of ADMM applied to quadratic programs; see Appendix~\ref{apdx:prox-general}.

\paragraph{Vectorization.}
If type group $\ell$ has a simple cost function whose proximal operator has an analytical expression, then it can be evaluated in parallel across all $|\ell|$ devices, the $T$ timesteps, and (for $\ell = \ell_1$) the $K+1$ contingencies.
For more complex device cost functions with proximal operators that must be evaluated using an approximate iterative algorithm, we can still parallelize the execution of the algorithm across the $|\ell|$ devices and the $K+1$ contingencies.

We implement this parallelism using vectorized operations.
For example, if the contingency devices are AC lines, we apply broadcasted operations to tensors of size $(K+1) \times |\ell_1| \times T$.
For devices like batteries, whose proximal operators are evaluated using a few steps of an iterative algorithm, we apply broadcasted operations to matrices of size $|\ell| \times T$ to evaluate all the battery proximal operators in parallel.
In a GPU computing language like PyTorch, these broadcasted operations are implemented by optimized GPU kernels that fully exploit the GPU's vector parallelism.

\subsection{Averages and Residuals}

For type $\ell$, let $A_{\ell i} \in \reals^{N \times |\ell|}$ be the matrix with entries,
\begin{equation*}
    (A_{\ell i})_{n d} = \begin{cases}
        1 & \textrm{terminal $i$ of device $d$ connects to node $n$}, \\
        0 & \textrm{otherwise}.
    \end{cases}
\end{equation*}
where $i = 1, \ldots, \tau_\ell$.
The matrices $A_{\ell i}$ are the \textit{incidence matrices} for each terminal of each device type.
Then averages such as $\bar p$ are given by
\begin{equation*}
    \bar p = \frac{1}{|n|} \odot \sum_{\ell} \sum_{i = 1}^{\tau_\ell} A_{\ell i} p_{\ell i},
\end{equation*}
where $p_{\ell i}$ is a $(K+1) \times |\ell| \times T$ tensor, $\bar p$ is a $(K+1) \times N \times T$ tensor, $|n| \in \reals^N$ is the vector of number of terminals per node, and $/$ and $\odot$ are elementwise division and multiplication, respectively.
Here, matrix-tensor multiplication is to be understood as matrix-matrix multiplication between the matrix and the last two dimensions of the tensor, broadcasted across all the previous dimensions of the tensor.
In our example, this means multiplying a single matrix of size $N \times |\ell|$ by $K+1$ matrices of size $|\ell| \times T$ to produce a tensor of size $(K+1) \times N \times T$.

Residuals such as $\tilde p$ are similarly given by,
\begin{equation*}
    \tilde p_{\ell i} = p_{\ell i} - A_{\ell i}^T \bar p.
\end{equation*}
Informally, averaging, and other operations that gather information by node, can be thought of as multiplication by a sparse incidence matrix.
On the other hand, computing residuals, and other operations that send node information to individual devices, can be thought of as multiplication by the adjoint of the same matrix.

On the GPU, multiplying by $A_{\ell i}$ can be implemented using a \textit{scatter} kernel.
Multiplying by $A_{\ell i}^T$ can be implemented with a \textit{gather} kernel.
Both of these operations are readily available and highly optimized in all major GPU computing languages, allowing us to efficiently implement average and residual operations on the GPU.

\subsection{Unrolled Differentiation}
\label{sec:unroll}

Message passing can be written as a sequence of operations,
\begin{equation*}
    s^\iter(x) = (h^\iter(\cdot, x) \circ h^\iterminus(\cdot, x) \circ \cdots \circ h^{(1)}(\cdot, x)) (s^{(0)}),
\end{equation*}
where $s^\iter(x) = (p^\iter, \theta^\iter, u^\iter, v^\iter)$, the function $h^\iter$ is the $i$th step of the message passing algorithm, and $x = (x^0, \ldots, x^K)$ is the problem data.
Here, $\circ$ is standard function composition.
Each iteration $h^\iter$ is a differentiable function, since the average and residual operators are linear and the proximal operator of a function in the form of~\eqref{eq:prox-general-form} is piecewise linear-quadratic on a polyhedral domain and, therefore, almost everywhere (a.e.) differentiable~\cite[Chapter 12C, Proposition 12.30]{rockafellar_variational_2009}.
In particular, the proximal operators for devices with simple cost functions are implemented as a composition of a.e. differentiable operations.
For more complex devices, proximal operators are implemented as a few iterations of a.e.\ differentiable operations.
Therefore, the entire algorithm is a.e.\ differentiable, since it consists of a composition of many a.e.\ differentiable functions, and we can compute the Jacobian $\partial s^\iter(x)$, the gradient of the state at iteration $i$ with respect to the problem data.

In a GPU computing language with \textit{automatic differentiation} support, we can enable \textit{gradient tracking} to automatically compute $\partial s^\iter(x)$.
Gradient tracking builds a dynamic computation graph and records intermediary values to automatically evaluate gradients.
In practice, the full Jacobian matrix is never instantiated explicitly, but rather used as an operator in \textit{reverse mode} or \textit{forward mode} differentiation~\citep{baydin2018automatic}.

\subsection{Software Implementation}

We implement our message passing algorithm in pure PyTorch.
The entire algorithm can be run on the GPU, leverage PyTorch's automatic differentiation capabilities, and seamlessly integrate into a larger PyTorch model, e.g., a physics-informed neural network~\citep{karniadakis2021physics, chen2021enforcing, amos2017optnet}.
Our software package is also modular: the user can specify new devices by writing code for that device's cost function and its proximal operator.
We provide code for devices such as generators, loads, AC and DC transmission lines, and batteries.
Complex proximal operators, such as for batteries, are evaluated by running 10 iterations of the ADMM-based quadratic program solver in Appendix~\ref{apdx:prox-general}.
On a single Nvidia A100 GPU, our code takes less than 4 seconds to run 1000 iterations of message passing for a network with 500 nodes, 3200 devices, 24 time periods, and 100 contingencies.
Similarly, when computing gradients using automatic differentiation, the backwards pass takes just under 4 seconds;
in most cases, we find that evaluating gradients takes a similar time to running the algorithm forward.
Our implementation is open-source and available on Github at \url{https://github.com/degleris1/zap}.

\section{Numerical Experiments}

In this section, we demonstrate the convergence properties and computational performance of messaging passing on contingency-constrained OPF problems.
We show empirically that our solver significantly outperforms a CPU-based commercial solvers, Mosek~\citep{mosek}, and explain how unrolled differentiation can be applied to solving expansion planning problems.

\subsection{Problem Data}

We use the open-source PyPSA-USA dataset~\citep{tehranchi2024pypsa} to generate realistic models of the U.S.\ Western Interconnect.
The network contains 4786 nodes and 8760 hours of data from 2019, with the option to change the network resolution to as few as 30 nodes.
The devices in the network include loads, generators, DC and AC transmission lines, and batteries.
We make several small adjustments to the dataset to create realistic grid cases;
however, the focus of our experiments is on computational performance and convergence rate, which should generally be unaffected by these details.
Further details on the dataset are available in Appendix~\ref{apdx:data}.

\subsection{Experiment Settings}

In all our experiments, we use the AC lines as the contingency devices and define each contingency consist of a single line outage.
we set the penalty parameters $\rho_p = \rho_\theta = 1$ and use the adaptive penalty parameter update from Section~\ref{sec:adaptive-penalty} with $\mu = 2$ and $\tau = 1.1$ every 10 iterations.
For simplicity, we ignore over-relaxation and set $\alpha = 1$.
We run all our experiments on a single A100 GPU with 80~Gb of memory supported by 32 virtual CPU cores (16 physical cores) and 64~Gb of RAM.

\subsection{Convergence}
\label{sec:exp-converge}

\begin{figure}[t]
    \centering
    \includegraphics[width=\fullFigWidth]{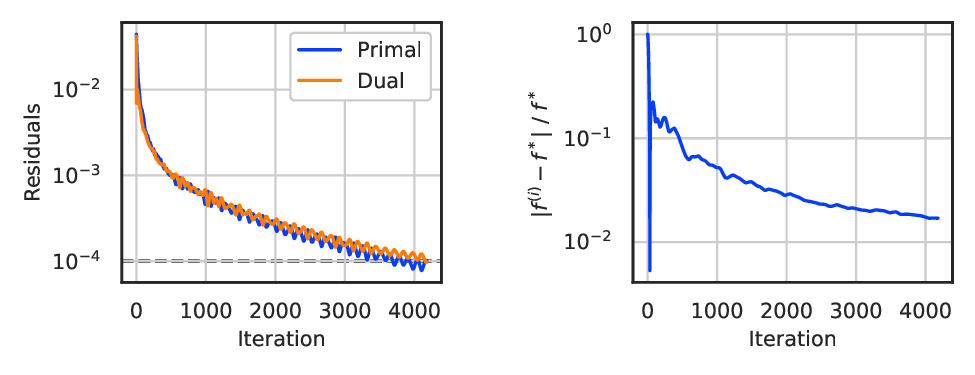}
    \caption{
        Primal and dual residuals (left) and relative objective objective suboptimality (right) as a function of algorithm iteration.
        Message passing converges to an accuracy of $10^{-3}$ in 529 iterations and $10^{-4}$ in 4180 iterations.
    }
    \label{fig:convergence}
\end{figure}

In our first experiment, we consider a simple problem with 500 nodes, 3112 devices, 24 hours, and no contingencies.
The devices include 500 loads, 1392 generators, 1143 AC transmission lines, 3 HVDC transmission lines, and 74 batteries.
We set the tolerance to $\epsilon = 10^{-4}$ and run the algorithm until convergence as defined in Equation~\eqref{eq:converge}.
In Figure~\ref{fig:convergence}, we plot the primal residuals, dual residuals, and relative objective suboptimality (to compute the optimal objective value, we model the same problem in CVXPY~\citep{diamond2016cvxpy} and solve it using Mosek).

Message passing takes just over 4000 iterations to converge to the specified tolerance.
However, the algorithm reaches an accuracy of $10^{-3}$ in just 529 iterations.
At convergence, message passing finds a solution with objective value $f^\iter = 181.4$, just 1.6\% different from the true optimal value $f^* = 184.5$.

\begin{figure}
    \centering
    \includegraphics[width=\smallFigWidth]{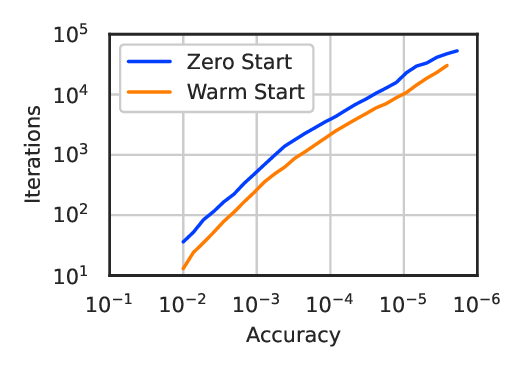}
    \caption{
        Number of iterations to converge as a function of tolerance $\epsilon$ when initializing with all zeroes (blue) and the solution to the previous $T$ time periods (orange).
        Warm starting reduces the required iterations by about a half, e.g., 279 iterations versus 514 iterations for an accuracy of $10^{-3}$.
    }
    \label{fig:accuracy}
\end{figure}

To better understand the rate of convergence for message passing, we run the same example for 100,000 iterations and plot the number of iterations until convergence as a function of accuracy in Figure~\ref{fig:accuracy}.
This is plotted for two different initializations: a zero start initialization where $(p^{(0)}, \theta^{(0)}, u^{(0)}, v^{(0)}) = 0$ and a warm start initialization where the solution to the previous 24 hours is used as the initial point.
For both initializations, we find that the resulting accuracy is inversely proportional to the number of iterations: to increase the accuracy by a factor of 10, we need to run roughly 10 times more iterations.
This suggests that message passing is best suited to finding low to medium accuracy solutions (reaching accuracies of, say, $10^{-2}$ to $10^{-4}$) and to applications where only an approximate solution is needed.

When initializing the algorithm using the previous day's solution, message passing converges in fewer iterations in all cases.
In particular, warm starting the solver in this manner generally reduces the number of iterations to converge to a particular accuracy by about 50\%.

\subsection{Computational Scaling}

\begin{figure}
    \centering
    \includegraphics[width=\fullFigWidth]{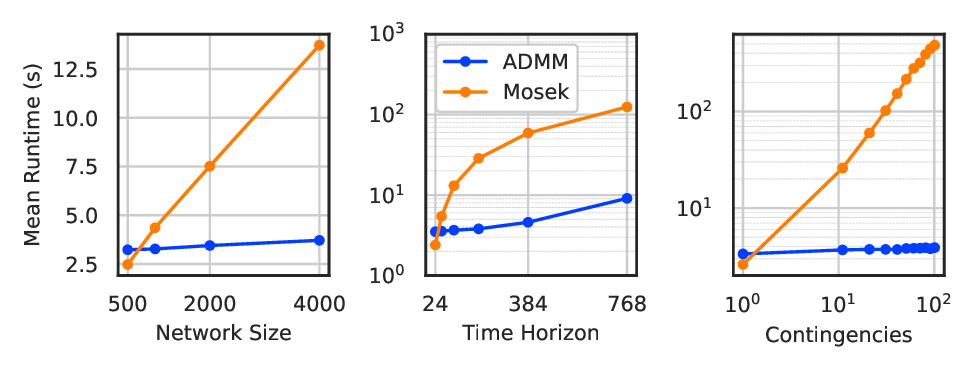}
    \caption{
        Runtimes for 1000 iterations of message passing compared with a commercial interior point method, Mosek, solved to a low accuracy of $10^{-3}$.
        Problem size is varied by increasing the number of nodes $N$ (left plot), the time horizon $T$ (middle plot), and the number of contingencies $K$ (right plot).
        In all settings, message passing is significantly faster than Mosek for sufficiently large cases, achieving over a 400$\times$ speedup in the largest of cases.
    }
    \label{fig:scaling}
\end{figure}

In this experiment, we study the computational scaling of message passing for different problem sizes.
In particular, we use the 500 node, 24 hour, no contingency case from Section~\ref{sec:exp-converge} as a starting point, then increase either the number of nodes (and thus devices) in the network, the time horizon, or the number of contingencies.
We record the runtime for 1000 iterations of message passing and compare it to Mosek, which is set with a relatively low accuracy of $10^{-3}$ for fair comparison.
All problem sizes are solved for three different cases: a low load case, a medium load case, and a high load case that involves load curtailment.
Runtimes are averaged across cases.
These results are shown in Figure~\ref{fig:scaling}.

In all settings, message passing is significantly faster than Mosek for sufficiently large cases.
In particular, message passing is over 400x faster than Mosek for the case with 500 contingencies, suggesting that our algorithm is particularly well suited to contingency-constrained OPF problems.
Note that in the 1000 contingency case, Mosek did not return a solution within a 2 hour time limit, so no point is generated;
in contrast, message passing takes about 20 seconds.

We track runtimes for a fixed number of iterations since the convergence rate can vary significantly based on problem specific details.
However, in order to ensure that message passing converges to a reasonable solution, we track the maximum of the root mean square (RMS) primal and dual residuals for each case.
The median RMS residual across all cases is $0.904 \cdot 10^{-3}$, and the largest is $3.035 \cdot 10^{-3}$.
For reference, the average hourly load across all cases is $67.44~\textrm{GW}$.
This means that the average node has a power imbalance of less than 1~MW in a problem with 67.44~GW of load.

\paragraph{Large example.}
To demonstrate the full potential of our method, we solve a large case with $N=4000$ nodes, $D = 17426$ devices, $T=24$ time periods, and $K = 1000$ contingencies.
The full problem has over 500 million variables in total.
For this example, running 1000 iterations of message passing takes about 60 seconds.
This demonstrates that our method is capable of quickly solving real world contingency-constrained OPF problems across multiple time periods coupled by ramping or energy storage constraints.
For even larger cases (e.g., $K = 10^{4}$), message passing exceeds the 80~Gb memory limit on the GPU;
this can be resolved by using high memory or potentially multiple GPUs.
Using multiple GPUs may involve expensive communication steps that require further analysis.
We leave this as future work discussed in Section~\ref{sec:conclusion}.

\subsection{Sensitivities}

\begin{figure}
    \centering
    \includegraphics[width=\fullFigWidth]{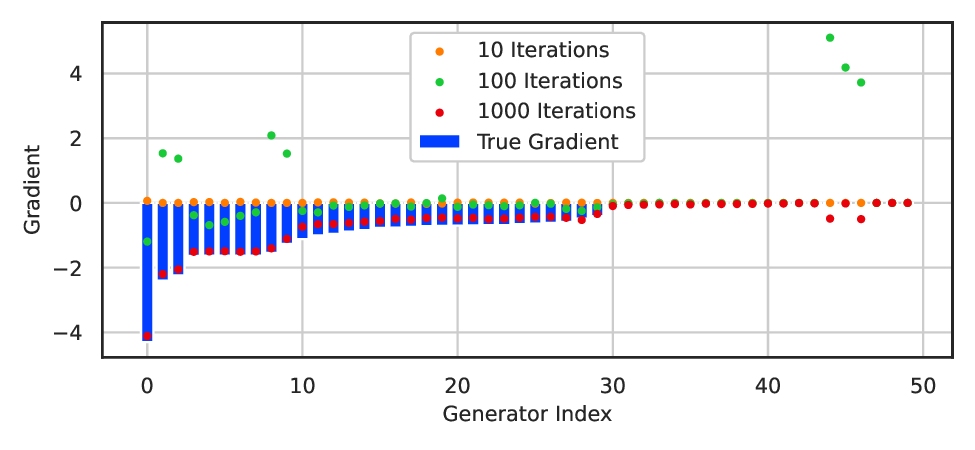}
    \caption{
        Derivatives of total cost with respect to generator capacities for 50 randomly selected generators.
        (Blue bars) Exact derivatives computed using Mosek and duality theory.
        (Orange, green, and red dots) Approximate derivatives from unrolling message passing with automatic differentiation for $10$, $100$, and $1000$ iterations.
        Unrolling for more iterations leads to more accurate derivatives.
    }
    \label{fig:derivatives}
\end{figure}

We track gradients for $10$, $100$, and $1000$ iterations of message passing using automatic differentiation.
We then compute $f^\iter = f(p^\iter, \theta^\iter; x)$ and evaluate the gradient of the total system cost with respect to generator capacities, $\partial f^\iter / \partial p^{\max} \in \reals^r$.
Since $\partial f^\iter / \partial p^{\max} = - \lambda$~\citep{boyd2004convex}, where $\lambda \in \reals^r_{+}$ is the dual variable of the generator capacity constraint, we can evaluate the derivatives exactly by solving OPF using Mosek and extracting the dual variables.
We use this as a baseline to measure the accuracy of the approximate derivatives produced by unrolling message passing.
We run this experiment for the problem in Section~\ref{sec:exp-converge} and plot the resulting sensitivities in Figure~\ref{fig:derivatives}.
As expected, running the algorithm for more iterations leads to a more accurate approximation of the derivatives.

\subsection{Application: Expansion Planning}

In our last experiment, we highlight how our algorithm can be applied to expansion planning problems.
We consider an expansion planning problem on the 500 node system over $M = 8$ scenarios and $T = 24$ hours.
We use the multi-value expansion planning framework from~\citet{degleris2024gradient} and solve a carbon-aware planning problem~\citep{degleris2021emissions} with no carbon tax but with an implicit carbon penalty of \$200 per mTCO$_2$.
The resulting optimization problem can be formulated as
\begin{equation}
\label{eq:planning}
\begin{array}{ll}
    \textrm{minimize} \quad
    & \sum_{m=1}^M F_m(s_m^*(\eta)) \\
    \textrm{subject to} \quad
    & \eta^{\min} \leq \eta \leq \eta^{\max},
\end{array}
\end{equation}
where the variable is the vector of device capacities $\eta \in \reals^R$.
Here, $s_m^*(\eta) = (p^*(\eta), \\ \theta^*(\eta), u^*(\eta), v^*(\eta))$ is the \textit{solution map} of scenario $m$, which gives the solution to Problem~\eqref{eq:admm-problem-full} as a function of the device capacities, and $F_m : \reals^{4JT(K+1)} \rightarrow \reals$ is the \textit{planner's objective}.
In our example, the device capacities are for all expandable AC lines, DC lines, generators, and batteries.
The planner's objective is the total system cost plus an \textit{implicit} penalty on carbon emissions;
the penalty is implicit because only the planner accounts for it in their objective, and the OPF problem itself has no explicit carbon tax~\citep{degleris2021emissions}.

Problem~\eqref{eq:planning} is a bilevel optimization problem, which is non-convex and NP-hard to solve in general.
We use gradient descent as an approximate local solution method~\citep{degleris2024gradient}: at each iteration, we evaluate $s^*(\eta)$ by solving~\eqref{eq:admm-problem-full} using message passing, then compute the gradient $\Delta = \sum_m \partial s_m^*(\eta)^T \cdot \nabla F_m(s_m^*(\eta))$ using unrolled differentiation, then update the capacities as
$
   \eta := \eta - \alpha \Delta,
$
where $\alpha$ is the step size.
In our experiments, we set $\alpha = 0.01$.
At each gradient step, we set our message passing algorithm to run for a minimum of 250 iterations and a maximum of 1000 iterations, using a tolerance $\epsilon = 10^{-3}$ and warm starting from the solution from the previous step.
In general, because the capacities $\eta$ only change slightly between gradient steps, message passing should only need a few iterations to find a good solution.
Finally, we implement a batch solver that can solve multiple OPF problems in parallel, so all 8 scenarios are solved in parallel on the same GPU.

\begin{figure}
\centering
\includegraphics[width=\fullFigWidth]{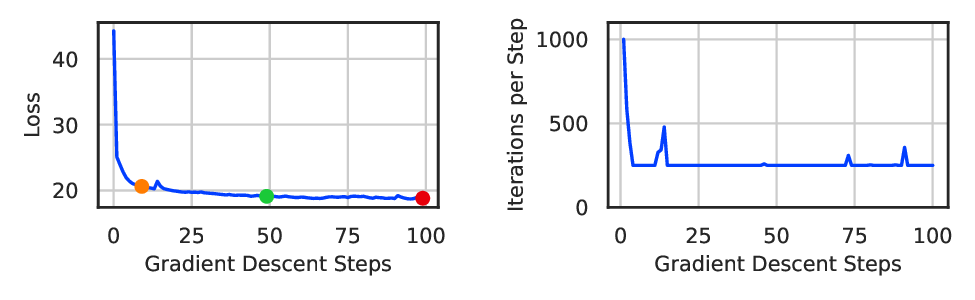}
\caption{
    Problem loss (left) and number of message passing iterations per gradient step (right) as a function of gradient descent step.
    Colored dots are problem losses at 10, 50, and 100 iterations (see Figure~\ref{fig:planning-result}).
    The message passing algorithm is set to run for at least 250 and at most 1000 iterations, stopping early if it reaches an accuracy of $\epsilon = 10^{-3}$;
    each solve is initialized using the solution at the previous step.
    After a few iterations, message passing nearly always converges in the minimum number of iterations.
    Gradient descent rapidly converges in a few dozen steps.
}
\label{fig:planning-metrics}
\end{figure}

\begin{figure}
\centering
\includegraphics[width=\fullFigWidth]{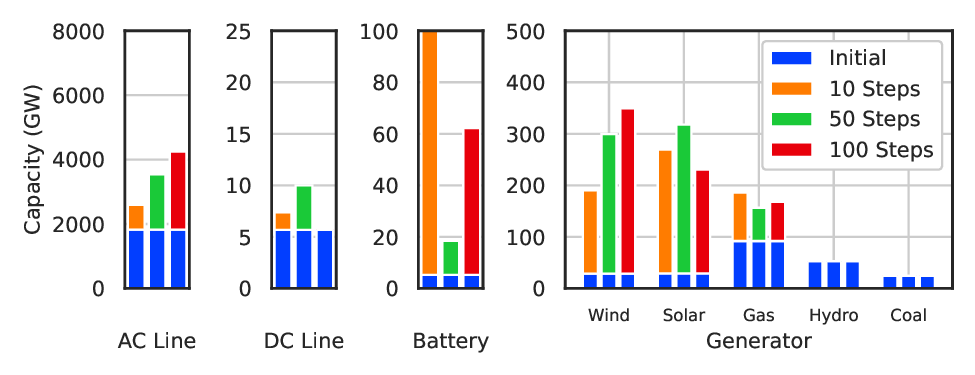}
\caption{
    Expansion in transmission line, battery, and generator capacities (hydro and coal expansion is not permitted) after 10, 50, and 100 iterations of gradient descent.
    An implicit carbon penalty of \$200 per mTCO$_2$ leads to a significant expansion of renewable generation and storage.
}
\label{fig:planning-result}
\end{figure}

We run gradient descent for 100 iterations and plot the problem loss and the number of message passing iterations per gradeint step in Figure~\ref{fig:planning-metrics}.
As anticipated, message passing usually runs for the minimum number of iterations after the first few steps.
Gradient descent rapidly converges in the first few tens of iterations, then slowly improves the solution thereafter.
We show the final expansion in device capacities in Figure~\ref{fig:planning-result}.

Gradient descent takes 4 minutes and 24 seconds to run in total.
For comparison, we ran the same gradient descent algorithm except using Mosek to solve the OPF problems and implicit differentiation to compute derivatives~\citep{degleris2024gradient}.
The two methods achieve nearly identical results.
However, the Mosek-based algorithm takes 1 hour and 26 minutes to run the same 100 gradient steps.

\section{Conclusion}
\label{sec:conclusion}

In this work, we present a message passing algorithm for contingency-constrained OPF problems.
Our algorithm leverages the sparsity inherent in electricity networks, does not require solving any linear systems, and runs entirely on the GPU.
We implement this algorithm in pure PyTorch and show that it performs well in finding low to medium accuracy solutions.
In particular, we achieve over 100$\times$ speedups compared to an off-the-shelf commerical solver, Mosek, solving problems with thousands of nodes and contingencies in under a minute.
Moreover, our algorithm consists only of vectorized proximal operators and sparse matrix multiplies and is therefore end-to-end differentiable: every iteration can be unrolled using automatic differentiation.
We demonstrate one potential application of our method to grid expansion planning.
In general, we believe our method could be applied to numerous settings in which OPF problems need to be solved to a modest accuracy at scale, such as stochastic production cost modeling, long-term expansion planning, and physics-informed machine learning.

Future work could consider applying GPU-accelerated message passing to the aforementioned applications, extending the algorithm to devices with non-convex device costs, or supporting multiple GPUs.
Non-convex device costs could be addressed, for example, by using sequential linearization, applying convex relaxations, or even evaluating the proximal operators of the non-convex device cost functions directly (when possible).
Supporting multiple GPUs would require analyzing the communication steps between GPUs, which is critical to understanding the algorithm's performance on decentralized architectures.

\appendix
\section{Notation}
\label{apdx:notation}

Because the contingency-constrained OPF involves a number of different variables and indices, we keep a comprehensive list of the notation here.
In general, the symbols listed here are exclusively used to denote a single type of quantity, e.g., $j$ is always a terminal and never an index for some other quantity.
We use the following symbols for sets and indices:
\begin{itemize}
    \item $j \in \mathcal J = \{1, \ldots, J\}$ for \textit{terminals},
    \item $n \in \mathcal N = \{ n_1, \ldots, n_N \}$ for \textit{nodes}, where $n \subset \mathcal J$ and $\mathcal N$ partitions $\mathcal J$,
    \item $d \in \mathcal D = \{ d_1, \ldots, d_D \}$ for \textit{devices}, where $d \subset \mathcal J$ and $\mathcal D$ partitions $\mathcal J$,
    \item $\ell \in \mathcal L = \{\ell_1, \ldots, \ell_L \}$ for \textit{type groups}, where $\ell \subset D$ and $\mathcal L$ partitions $\mathcal D$,
    \item $k \in \{ 0, \ldots, K\}$ for contingencies,
\end{itemize}
We use the following symbols for variables:
\begin{itemize}
    \item $p \in \reals^{(K+1) \times J \times T}$ for power (injection) schedules,
    \item $\theta \in \reals^{(K+1) \times J \times T}$ for phase (angle) schedules,
    \item $u \in \reals^{(K+1) \times J \times T}$ for the scaled dual variables of the power balance constraints,
    \item $v \in \reals^{(K+1) \times J \times T}$ for the scaled dual variables of the phase consistency constraints.
\end{itemize}
We use set-valued indices flexibly to denote submatrices or subtensors.
\begin{itemize}
    \item $p_{kj} \in \reals^T$ refers to the contingency-$k$ power schedule for a single terminal during contingency $k$,
    \item $p_{kn} \in \reals^{|n| \times T}$ refers to the contingency-$k$ power schedule for all the terminals connected to node $n$,
    \item $p_{kd} \in \reals^{|n| \times T}$ refers to the contingency-$k$  power schedule for all the terminals connected to device $d$,
    \item $p_{k \ell} \in \reals^{|\ell| \times |\tau_\ell| \times T}$ refers to the contingency-$k$ power schedule for all the terminals connected to a device in type group $\ell$.
\end{itemize}
We drop the first index, the contingency, when there are no contingencies involved. 
We use the following for device data and algorithm parameters:
\begin{itemize}
    \item $x_\ell \in \reals^{|\ell| \times r_\ell}$, where $\ell \in \{\ell_1, \ldots, \ell_L\}$, for problem data,
    \item $\rho_p \in \reals_+$ and $\rho_\theta \in \reals_+$ for the penalty parameters for the power balance and phase consistency constraints, respectively.
\end{itemize}
Each device has a cost function $f_d(p_d, \theta_d)$ that is defined by its device group and problem data;
a device $d \in \ell$ with data $x_d \in \reals^{r_\ell}$ has cost function
\begin{equation*}
    f_d(p_d, \theta_d) = g_\ell(p_d, \theta_d; x_d).
\end{equation*}
We also refer to the sum of device costs in a type group, $G_\ell(p_\ell, \theta_\ell) = \sum_{d \in \ell} g_\ell(p_d, \theta_d; x_d)$.

\section{Averages and Residuals}
\label{apdx:average-resid}

\newcommand{\average}{\mathbf{avg}}
\newcommand{\resid}{\mathbf{resid}}

Here we state basic properties of the node average operator,
\begin{equation*}
    \bar z_j = \average(z_j) = \frac{1}{|n|} \sum_{j' \in n} z_{j'},
\end{equation*}
where $j$ is a terminal connected to node $n$, and the node residual operator,
\begin{equation*}
    \tilde z_j = \resid(z_j) = z_j - \bar z_j.
\end{equation*}
First, both the average and residual operators are linear.
In particular, the average operator can be implemented with the matrix $A = L^T \diag(|n|)^{-1} L$ where $L \in \reals^{N \times J}$ is the node-terminal incidence matrix and $\diag(|n|) \in \reals^{N \times N}$ is a diagonal matrix whose entries are the number of terminals connected to each node.
Likewise, the residual operator is given by the matrix $R = I - A$.
Both matrices are symmetric, $A^T = A$ and $R^T = R$.

Second, both the average and residual operators are idempotent,
\begin{align*}
    \average(\average(z_j)) &= \average(z_j), \\
    \resid(\resid(z_j)) &= \resid(z_j).
\end{align*}
Third, the average and residual operators cancel one another,
\begin{align*}
    \average(\resid(z_j)) = \resid(\average(z_j)) = 0.
\end{align*}
Finally, we note that
\begin{align*}
    \resid(z_j) = 0 &\quad\longleftrightarrow\quad z_j = \average(z_j), \\
    \average(z_j) = 0 &\quad\longleftrightarrow\quad z_j = \resid(z_j).
\end{align*}
All the above properties can be checked with simple algebra.

\section{Dual Problem}
\label{apdx:dual}

We consider the case when $K = 0$ and $T=1$ for simplicity.
Recall that Problem~\eqref{eq:dispatch} is
\begin{equation*}
\begin{array}{lll}
    \textrm{minimize} \quad
    & \sum_d f_d(p_d, \theta_d)
    \\[1em]
    \textrm{subject to} \quad
    & \bar p = 0, \\
    & \tilde \theta = 0, \\
\end{array}
\end{equation*}
where the variables are $p, \theta \in \reals^{J}$.
The Lagrangian is
\begin{align*}
    L(p, \theta, u, v) 
    &= u^T \bar p + v^T \tilde \theta + \sum_d f_d(p_d, \theta_d) \\
    &= \sum_d \left( f_d(p_d, \theta_d) 
        + \sum_{j \in d} \left( 
        v_j \theta_j - \sum_{j' \in n_j} \frac{1}{|n_j|} v_{j'} \theta_j
        + \sum_{j' \in n_j} \frac{1}{|n_j|} u_{j'} p_j
        \right) \right) \\
    &= \sum_d \left( f_d(p_d, \theta_d) 
        + \tilde v_j^T \theta_j
        + \bar u_j^T p_j
         \right),
\end{align*}
where $u, v \in \reals^{J}$. 
So the dual function is
\begin{align*}
    g(u, v) 
    = \inf_{p, \theta} L(p, \theta, u, v)
    &= \sum_d  \inf_{p_d, \theta_d} \left( f_d(p_d, \theta_d) 
        + \tilde v_j^T \theta_j
        + \bar u_j^T p_j\right) \\
    &= \sum_d -f_d^*( -\bar u_d, -\tilde v_d).
\end{align*}
where $f_d^* : \reals^{|d|} \times \reals^{|d|} \rightarrow \reals$ is the convex conjugate of $f_d$.
So the dual problem is
\begin{equation}
\label{eq:dual-simple}
\begin{array}{lll}
    \textrm{maximize} \quad
    & \sum_d f_d^*(-\bar u_d, -\tilde v_d), \\[0.5em]
\end{array}
\end{equation}
with variables $u, v \in \reals^J$.
Note that the objective is only function of $\bar u$ and $\tilde v$.
Moreover, for any vector $\bar u$, there is another vector $u'$ such that $\bar u = \bar u' = u'$ (of course, this vector is just $\bar u$ itself).
So the problem remains unchanged if we add the constraint $u = \bar u$ or, equivalently, $\tilde u = 0$.
Using the same logic, we can add the constraint $v = \tilde v$ or, equivalently, $\bar v = 0$.
Therefore, Problem~\eqref{eq:dual-simple} is equivalent to
\begin{equation}
\label{eq:dual-problem}
\begin{array}{lll}
    \textrm{maximize} \quad
    & \sum_d f_d^*(-u_d, -v_d)
    \\[1em]
    \textrm{subject to} \quad
    & \tilde u = 0, \\
    & \bar v = 0,
\end{array}
\end{equation}
with variables $u, v \in \reals^J$.
We make two remarks.
First, the dual of the DC optimal power flow problem is itself a DC optimal power flow problem, where the device costs have been replaced by their convex conjugates, the power balance constraint is now on the dual variables of the phase consistency constraints, and the phase consistency constraint is now on the dual variables of the power balance constraints.

Second, any solution to the dual problem~\eqref{eq:dual-problem} must satisfy $\tilde u^* = 0$ and $\bar v^* = 0$.
However, in some sense, these constraints are arbitrary: if $(u^*, v^*)$ is a solution to~\eqref{eq:dual-problem}, then any $(u, v)$ such that $\bar u = \bar u^*$ and $\tilde v = \tilde v^*$ is a solution to~\eqref{eq:dual-simple}.
This ultimately reflects the redundancy in the power balance and phase consistency constraints.
The $J$ power balance constraints are of rank $N < J$, and the $J$ phase consistency constraints are of rank $J - N$.
Therefore, we add the constraints $\tilde u = 0$ and $\bar v = 0$ at no loss of generality; problems~\eqref{eq:dual-simple} and~\eqref{eq:dual-problem} are equivalent.

\section{Evaluating Proximal Operators}
\label{apdx:prox-general}

Many device cost functions have proximal operators with simple closed form expressions.
In general, however, we consider device cost functions of the form
\begin{align}
\label{eq:prox-general-form-apdx}
    f(p, \theta) = \min_{s} \left( s^T Q s + q^T s + \indicator\{ A_1 p + A_2 \theta + A_3 s \leq b \} \right),
\end{align}
which has proximal operator
\begin{align}
\label{eq:prox-op-general}
    \prox_{f, \rho}(z, \xi) 
    &= 
    \argmin_{p, \theta} \left( 
        \min_{s} \left( s^T Q s + q^T s + \indicator\{ A_1 p + A_2 \theta + A_3 s \leq b \} \right) 
        + \frac{\rho_p}{2} \| p - z \|_2^2
        + \frac{\rho_\theta}{2} \| \theta - \xi \|_2^2
    \right),
\end{align}
which does not have a closed form solution.
To evaluate this operator, we use a classic form of ADMM for solving quadratic programs.
We introduce variables $\alpha, \beta \in \reals^{m}$ and solve the equivalent problem
\begin{equation*}
\label{eq:prox-admm-problem}
\begin{array}{ll}
    \textrm{minimize} \quad
    & s^T Q s + q^T s 
    + \frac{\rho_p}{2} \| p - z \|_2^2
    + \frac{\rho_\theta}{2} \| \theta - \xi \|_2^2 \\
    & \quad + \indicator\{ A_1 p + A_2 \theta + A_3 s - b - \alpha = 0 \}
        + \indicator \{ \beta \leq 0 \} \\[0.5em]
    \textrm{subject to} \quad 
    & \alpha = \beta,
\end{array}
\end{equation*}
with variables $p, \theta \in \reals^T$, $s \in \reals^\mu$, and $\alpha, \beta \in \reals^{m}$.
Then, the ADMM updates for~\eqref{eq:prox-admm-problem} are
\newcommand{\iterkplus}{{(k+1)}}
\newcommand{\iterk}{{(k)}}
\begin{enumerate}
    \item \textit{$(p, \xi, s, \alpha)$-update}:
    \begin{align*}
        (p^\iterkplus, \xi^\iterkplus, s^\iterkplus, \alpha^\iterkplus)
        =
        \argmin_{p, \xi, s, \alpha} \Big( 
            &s^T Q s + q^T s 
            + \frac{\rho_p}{2} \| p - z \|_2^2
            + \frac{\rho_\theta}{2} \| \theta - \xi \|_2^2 \\
            & + \indicator\{ A_1 p + A_2 \theta + A_3 s - b - \alpha = 0 \} \\
            &+ \frac{\omega}{2} \| \alpha - \beta^\iterk + \lambda^\iterk \|_2^2 
        \Big).
    \end{align*}
    \item \textit{$\beta$-update}:
    \begin{align*}
        \beta^\iterkplus
        &=
        \argmin_{\beta} \Big( 
            \indicator\{ \beta \leq 0 \}
            + \frac{\omega}{2} \| \beta - \alpha^\iterkplus - \lambda^\iterk \|_2^2.
        \Big) \\
        &= \min\left( \alpha^\iterkplus + \lambda^\iterk, 0 \right)
    \end{align*}
    \item \textit{Dual update:}
    \begin{align*}
        \lambda^\iterkplus =
        \lambda^\iterk + 
        \alpha^\iterk - \beta^\iterk.
    \end{align*}
\end{enumerate}
Here, $(k)$ is the iteration number and $\omega$ is the penalty parameter.
Updates~2 and~3 have closed form solutions that involve simple operations, and update~1 requires solving a quadratic program with linear equality constraints, which can be computed by solving a linear system of size $O(T+\mu+m)$~\citep{boyd2004convex}.
This linear system remains constant so long as the the problem data $(Q, q, A_1, A_2, A_3, b)$ and the penalty terms $\rho_p, \rho_\theta, \omega$ remain constant;
therefore, the factorization can be computed once and reused throughout the message passing algorithm.
In our implementation, we run this algorithm for 10 to 50 iterations to approximately evaluate~\eqref{eq:prox-general-form-apdx} during each message passing step.

\section{Simulation Data}
\label{apdx:data}

We use the open-source PyPSA-USA dataset~\citep{tehranchi2024pypsa} available at
\begin{equation*}
    \textrm{\url{https://github.com/pypsa/pypsa-usa}}
\end{equation*}
(we specifically use commit version v0.1.0 from April 24th, 2024).
We generate one year of hourly data (8760 hours) for the U.S.\ Western Interconnect.
Renewable generation and load is based on 2019 weather. 
Existing generator capacities comes from publicly available EIA generation data, and transmission data is based on the Breakthrough Energy network.
Load growth, fuel costs, and capital costs (for the expansion planning study) all use 2050 projections from the NREL ATB study.
We use the default settings to cluster the network to sizes between 500 and 4000 nodes.

After generating the data, we decrease generation and transmission capacities by 30\% to model retirements and outages.
Because load curves are for 2050, we scale load to between 50\% and 70\% of the projected value so that the cases can be solved without expansion. 
All load is curtailable with a curtailment cost of \$500 per MWh, and only the 70\% load scenario involves curtailment.

\section*{Acknowledgements}
Anthony Degleris is supported by the U.S. Department of Energy, Office of Science, Office of Advanced Scientific Computing Research, Department of Energy Computational Science Graduate Fellowship under Award Number DE-SC0021110.

This research used resources of the National Energy Research Scientific Computing Center (NERSC), a U.S. Department of Energy Office of Science User Facility located at Lawrence Berkeley National Laboratory, operated under Contract No. DE-AC02-05CH11231 using NERSC award DDR-ERCAP0026889.

\paragraph{}

\bibliography{main}

\end{document}